\documentclass[11pt]{amsart}
\usepackage{amssymb}

\theoremstyle{plain}
\newtheorem{theorem}{Theorem}
\newtheorem{corollary}[theorem]{Corollary}
\newtheorem{proposition}[theorem]{Proposition}
\newtheorem{lemma}[theorem]{Lemma}

\theoremstyle{remark}

\newcommand{\A}{\mathcal A}
\newcommand{\EA}{E(\mathcal A)}
\newcommand{\B}{\mathcal B}
\newcommand{\EB}{E(\mathcal B)}
\newcommand{\EX}{E(C(X))}
\newcommand{\EY}{E(C(Y))}
\renewcommand{\Re}{\operatorname{Re}}
\renewcommand{\Im}{\operatorname{Im}}
\newcommand{\half}{\frac{1}{2}}
\newcommand{\PA}{P(\mathcal A)}
\newcommand{\PB}{P(\mathcal B)}

\begin{document}
\title[]{Sequential isomorphisms between the sets of von Neumann algebra effects}
\author{LAJOS MOLN\'AR}
\address{Institute of Mathematics\\
         University of Debrecen\\
         4010 Debrecen, P.O.Box 12, Hungary}
\email{molnarl@math.klte.hu}
\urladdr{http://www.math.klte.hu/\~{}molnarl/}
\thanks{This paper was written when the author held
a Humboldt Research Fellowship. He is very grateful to the
Alexander von Humboldt Foundation for providing ideal conditions
for research and to his host Werner Timmermann (TU Dresden,
Germany) for very warm hospitality.
         The author also acknowledges support from
         the Hungarian National Foundation for Scientific Research
         (OTKA), Grant No. T043080, T031995, and from
         the Ministry of Education, Hungary, Grant
         No. FKFP 0349/2000.}
\date{September 5, 2003}
\subjclass{Primary: 46L60, 47B49}
\begin{abstract}
In this paper we describe the structure of all sequential
isomorphisms between the sets of von Neumann algebra effects. It
turns out that if the underlying algebras have no commutative
direct summands, then every sequential isomorphism between the
sets of their effects extends to the direct sum of a
*-isomorphism and a *-antiisomorphism between the underlying von
Neumann algebras.
\end{abstract}
\maketitle

\section{Introduction and Statement of the Result}

Effects are well-known to play important role in certain parts of
quantum mechanics, for instance, in the quantum theory of
measurement \cite{BusLahMit91}. In the Hilbert space formalism of
quantum mechanics the effects are the positive (bounded linear)
operators on a Hilbert space $H$ which are majorized by the
identity $I$. The set $E(H)$ of all (Hilbert space) effects on
$H$ is usually equipped with certain algebraic operations and/or
relations each of them having physical content and hence one
obtains different algebraic structures (cf. \cite[Section
V]{Lud83}).

In the present paper we are interested in the structure
corresponding to the so-called sequential product. This concept
was recently introduced by Gudder and Nagy in \cite{GudNag01}.
Roughly speaking, if $A,B$ are effects, then their sequential
product $A\circ B$ as an effect corresponds to the sequential
measurement in which the measurement corresponding to $A$ is
performed first and the one corresponding to $B$ is performed
second. Mathematically, the definition is $A\circ B=A^{1/2} B
A^{1/2}$ $(A,B\in E(H))$. Just as with any other algebraic
structure, it is an important problem to explore the
corresponding automorphisms (which are called sequential
automorphisms) and to describe their structure if possible. The
problem concerning the whole set $E(H)$ has been recently solved.
In fact, as an easy consequence of our result \cite[Theorem
3]{ML00f}, Gudder and Greechie proved in \cite[Theorem
1]{GudGre02} that, supposing $\dim H\geq 3$, the sequential
automorphisms of $E(H)$ are exactly the transformations $\phi$ of
the form
\[
\phi(A)=UAU^* \qquad (A\in E(H)),
\]
where $U$ is an either unitary or antiunitary operator on $H$. We
note that as a byproduct of one of our results concerning certain
preservers on Hilbert space effects, in \cite[Corollary 7]{MLxx7}
we obtained that the same result holds also when $\dim H=2$ (the
one-dimensional case is obviously exceptional).

It is a natural step in any research on operator algebras that if
one has a result concerning $B(H)$ (the algebra of all bounded
linear operators on $H$), then one is tempted to try to extend it
in some way for more general operator algebras. The first and in
many cases probably the most important candidates for such a
purpose are the von Neumann algebras. Therefore, after describing
the sequential automorphisms of the set of all Hilbert space
effects it is a natural problem to investigate the question for
von Neumann algebra effects, that is, for the set $\EA$ of all
elements of a given von Neumann algebra $\A$ which are positive
and majorized by the identity. (Observe that we have
$\EA=E(H)\cap \A$, where $H$ is the Hilbert space on which the
elements of $\A$ act.) We already have a result in this
direction. In fact, once again as a byproduct of an investigation
on a different problem, we obtained in \cite[Theorem 3]{MLxx5}
that if $\A$ is a von Neumann factor which is not of type I$_1$,
I$_2$, then every sequential automorphism of $\EA$ extends either
to a *-automorphism or to a *-antiautomorphism of $\A$.

In the present paper we go much further and give a complete
answer to the problem. Namely, we describe the structure of all
sequential isomorphisms between the sets of general von Neumann
algebra effects without making any restriction on the underlying
algebras.

We begin with the notation and some definitions that we shall use
throughout the paper. If $\A$ is a unital $C^*$-algebra, then the
effects in $\A$ are the positive elements of $\A$ which are less
than or equal to the unit of $\A$. The set of all effects in $\A$
is denoted by $E(\A)$. The sequential product $\circ$ on $\EA$ is
defined by
\[
A\circ B=A^{1/2} B A^{1/2} \qquad (A, B\in \EA).
\]
(For some interesting properties of this operation on $E(H)$
concerning associativity and commutativity see \cite{GudNag01}.)
If $\A, \B$ are unital $C^*$-algebras, then a bijective map
$\phi:\EA \to \EB$ which satisfies
\[
\phi(A\circ B)=\phi(A)\circ \phi(B) \qquad (A,B \in \EA)
\]
is called a sequential isomorphism.

Let $\psi_1 :\A \to \B$ be a *-isomorphism, that is, a bijective
linear map such that
\[
\psi_1(AB)=\psi_1(A)\psi_1(B), \quad \psi_1(A^*)=\psi_1(A)^*
\quad (A,B \in \A).
\]
It is easy to see that the restriction of $\psi_1$ onto $\EA$ is
a sequential isomorphism from $\EA$ onto $\EB$. Similar assertion
applies for any *-antiisomorphism $\psi_2: \A \to \B$ as well.
The bijective linear map $\psi_2: \A \to \B$ is called a
*-antiisomorphism if it satisfies
\[
\psi_2(AB)=\psi_2(B)\psi_2(A), \quad \psi_2(A^*)=\psi_2(A)^*
\quad (A,B \in \A).
\]
Clearly, the direct sum of sequential isomorphisms is again a
sequential isomorphism. Therefore, the direct sum of a
*-isomorphism and a *-antiisomorphism is, when restricted onto
the set of effects, a sequential isomorphism. In our main result
which follows we see that if the underlying algebras $\A, \B$ are
von Neumann algebras, then every sequential isomorphism between
the sets of their effects can be obtained in this way on the
'noncommutative parts' of $\EA$ and $\EB$.

\begin{theorem}\label{T:hetfo}
Let $\A,\B$ be von Neumann algebras and let $\phi:\EA \to \EB$ be
a sequential isomorphism. Then there are direct decompositions
\[
\A =\A_1 \oplus \A_2 \oplus \A_3 \quad \text{and} \quad
\B =\B_1 \oplus \B_2 \oplus \B_3
\]
within the category of von Neumann algebras and there are
bijective maps
\[
\phi_1:E(\mathcal A_1) \to E(\mathcal B_1), \quad
\Phi_2:\A_2 \to \B_2, \quad
\Phi_3:\A_3 \to \B_3
\]
such that
\begin{itemize}
\item[(i)] $\A_1, \B_1$ are commutative von Neumann algebras and
the algebras $\A_2 \oplus \A_3$, $\B_2 \oplus \B_3$ have no
commutative direct summands;

\item[(ii)] $\phi_1$ is a multiplicative bijection, $\Phi_2$ is a
*-isomorphism, $\Phi_3$ is a *-antiisomorphism and $\phi=\phi_1
\oplus \Phi_2 \oplus \Phi_3$ holds on $\EA$.
\end{itemize}
(Of course, the above decomposition is meant in the sense that
some of the direct summands can be missing.) If there is a scalar
$\lambda \in ]0,1[$ such that $\phi(\lambda I)=\lambda I$, or
$\A$ has no commutative direct summand, then $\phi$ extends to
the direct sum of a *-isomorphism and a *-antiisomorphism between
subalgebras of $\A$ and $\B$.
\end{theorem}

To see that between the 'commutative parts' of the underlying
algebras the behavior of $\phi$ can be quite 'irregular',
consider an arbitrary compact Hausdorff space $X$. Let $C(X)$
denote the algebra of all continuous complex valued functions on
$X$. Take any strictly positive continuous function $p:X \to [0,
\infty [$ and define $\phi: E(C(X)) \to E(C(X))$ by
\[
\phi(f)(x)={f(x)}^{p(x)} \qquad (x\in X, f\in E(C(X))).
\]
It is easy to verify that $\phi$ is a sequential automorphism of
$E(C(X))$ which, in general, does not extend to an automorphism
of $C(X)$. (Observe that if $\A, \B$ are commutative, then the
sequential product on $\EA,\EB$ coincides with the ordinary
product and hence in that case the sequential isomorphisms
between $\EA$ and $\EB$ are just the multiplicative bijections
or, in other words, the semigroup isomorphisms.)

As for the case of factors, our main result has the following
immediate corollary.

\begin{corollary}
Let $\A, \B\neq \mathbb C I$ be factors and let $\phi:\EA \to
\EB$ be a sequential isomorphism. Then $\phi$ extends either to a
*-isomorphism or to a *-antiisomorphism between the algebras $\A$
and $\B$.
\end{corollary}

\section{Proof}

We present some further notation and definitions that we shall
use in the proof of our main result. If $\A$ is a $C^*$-algebra,
then $\A_s, \A^+$ denote the set of all self-adjoint elements of
$\A$ and the set of all positive elements of $\A$, respectively.
The set of all projections (i.e., self-adjoint idempotents) in
$\A$ is denoted by $\PA$ and $\mathcal Z(\A)$ stands for the
center of $\A$.

Our first auxiliary result concerns the so-called
$E$-isomorphisms of the sets of $C^*$-algebra effects. Let
$\mathcal A, \mathcal B$ be $C^*$-algebras. Let $\phi: \EA \to
\EB$ be a bijective map with the properties that
\[
A+B \in \EA \Longleftrightarrow \phi(A)+\phi(B) \in \EB \qquad
(A,B\in \EA)
\]
and
\begin{equation}\label{E:het5}
\phi(A+B)=\phi(A)+\phi(B) \qquad (A,B, A+B\in \EA).
\end{equation}
Then $\phi$ is called an $E$-isomorphism (cf.
\cite{CasDeVLahLev}).

To formulate our result on $E$-isomorphisms we also need the
concept of Jordan *-isomorphisms. Let $\A, \B$ be *-algebras and
let $\phi :\A \to \B$ be a bijective linear map with the
properties that
\[
\phi(A^2)=\phi(A)^2, \quad \phi(A^*)=\phi(A)^* \quad (A\in \A).
\]
Then $\phi$ is called a Jordan *-isomorphism. It is easy to see
that in the above definition the condition that
$\phi(A^2)=\phi(A)^2$ holds for every $A\in \A$ can be replaced
by the following one:
\[
\phi(AB+BA)=\phi(A)\phi(B)+\phi(B)\phi(A) \qquad (A, B\in \A).
\]
One can readily verify that if $\A,\B$ are $C^*$-algebras and
$\phi: \A \to \B$ is a Jordan *-isomorphism, then its restriction
onto $\EA$ is an $E$-isomorphism from $\EA$ onto $\EB$. Our first
result says that the converse is also true, that is, every
$E$-isomorphism from $\EA$ onto $\EB$ extends to a Jordan
*-isomorphism from $\A$ onto $\B$.

\begin{proposition}\label{P:het1}
Let $\mathcal A, \mathcal B$ be $C^*$-algebras and let $\phi: \EA
\to \EB$ be an $E$-isomorphism. Then there is a Jordan
*-isomorphism $\Phi: \A \to \B$ such that
\[
\phi(A)=\Phi(A) \qquad (A\in \EA).
\]
\end{proposition}

\begin{proof}
The idea of the proof is easy. Using a quite elementary argument,
we can extend $\phi$ to a linear transformation $\Phi:\A \to \B$.
We next show that $\Phi$ is bijective, unital and preserves the
order in both directions. Finally, we apply a well-known result
of Kadison on order isomorphisms between $C^*$-algebras to
conclude that $\Phi$ is a Jordan *-isomorphism.

Turning to the details, first we show that $\phi$ preserves the
order. Let $A,B \in \EA$ be such that $A\leq B$. Then we have
$B=A+(B-A)$, where $A, B-A\in \EA$. As $\phi$ is an
$E$-isomorphism, we obtain
\[
\phi(B)=\phi(A)+\phi(B-A)\geq \phi(A).
\]
Therefore, $\phi$ preserves the order. Since $\phi^{-1}$ has
similar properties as $\phi$, it follows that $\phi$ preserves
the order in both directions. We easily deduce that $\phi(0)=0$
and $\phi(I)=I$.

Let $A\in \EA$ be arbitrary. By the partial additivity of $\phi$
(see \eqref{E:het5}) we have
\[
\phi(A)=\phi\biggl(n \frac{1}{n} A\biggr)=n
\phi\biggl(\frac{1}{n} A\biggr)
\]
which implies
\[
\phi\biggl(\frac{1}{n} A\biggr)=\frac{1}{n} \phi(A)
\]
for every $n\in \mathbb N$. If $n,k\in \mathbb N$ and $k\leq n$,
using the partial additivity of $\phi$ again, we obtain that
\[
\frac{k}{n} \phi(A)=k\phi\biggl(\frac{1}{n} A\biggr)=\phi\biggl(
\frac{k}{n} A\biggr).
\]
By the order preserving property and the just proved
rational-homogeneity of $\phi$ one can easily verify that
\[
\phi(\lambda A)=\lambda \phi(A)
\]
holds for every $A\in \EA$ and $\lambda \in [0,1]$.

Now it requires only elementary arguments to check that the
transformation $\psi_1 :\A^+ \to \B^+$ defined by
\[
\psi_1(A)=
\left\{%
\begin{array}{ll}
    \| A\| \phi\bigl(\frac{A}{\| A\|}\bigr), & \hbox{$0\neq A\in \A^+$;} \\
    0, & \hbox{$A=0$} \\
\end{array}%
\right.
\]
is bijective, additive, positive homogeneous, unital, preserves
the order in both directions and extends $\phi$. We omit the
details. To proceed, we define the transformation $\psi_2: \A_s
\to \B_s$ by
\[
\psi_2(A)=\psi_1(A^+)-\psi_1(A^-) \qquad (A\in \A_s)
\]
where $A^+,A^-$ denote the positive part and the negative part of
$A$, respectively. It is easy to verify that $\psi_2$ is
bijective, linear, unital, preserves the order in both directions
and extends $\psi_1$. Finally, let $\Phi :\A \to \B$ be defined
by
\[
\Phi(A)=\psi_2(\Re A)+i \psi_2(\Im A) \qquad (A\in \A)
\]
where $\Re A$ and $\Im A$ denote the real part and the imaginary
part of $A$, respectively. One can check that $\Phi$ is
bijective, linear, unital, preserves the order in both
directions, extends $\psi_2$ and hence extends $\phi$ as well. We
now refer to a result of Kadison (see \cite[10.5.32.
Exercise]{KadRin86}) stating that any bijective unital linear
transformation between $C^*$-algebras which preserves the order
in both directions is necessarily a Jordan *-isomorphism.
Applying this result to $\Phi$ we complete the proof.
\end{proof}

This simple result will play important role in the proof of our
main theorem. In fact, the basic idea of that proof is the
following. Let $\A, \B$ be von Neumann algebras and let $\phi:
\EA \to \EB$ be a sequential isomorphism. Consider the type
decomposition of $\A$ and $\B$ (see, for example, \cite[6.5.2.
Theorem]{KadRin86}). We shall see that $\phi$ maps the effects in
the type I$_n$ direct summand of $\A$ to the effects in the type
I$_n$ direct summand of $\B$. Therefore, roughly speaking, we can
consider our problem separately on type I$_n$-algebras and on
algebras without type I direct summands and then take direct
sums. It will turn out that on such algebras (with the exception
of type I$_1$ algebras), the sequential isomorphisms are all
$E$-isomorphisms and hence we can apply Proposition~\ref{P:het1}.
This is the plan of the proof. The details follow.

First we present a result concerning type I$_1$ (i.e.,
commutative) algebras. As we have already mentioned, on such
algebras the sequential isomorphisms coincide with the
multiplicative bijections, i.e., the semigroup isomorphisms.

\begin{proposition}\label{P:het2}
Let $X,Y$ be compact Hausdorff spaces. Let $\phi: \EX \to \EY$ be
a bijective multiplicative map. Suppose that there is a constant
$\lambda \in ]0,1[$ such that $\phi(\lambda)=\lambda$. Then
$\phi$ extends to a *-isomorphism $\Phi: C(X) \to C(Y)$.
\end{proposition}

\begin{proof}
The basic idea of the proof that one should transfer
multiplicativity to additivity is due to P. \v Semrl.

Let $g,f\in \EX$, $g\leq f$ and suppose that $f$ is nowhere
vanishing. Then $\frac{g}{f}\in \EX$, so there is a function
$f'\in \EX$ such that $g=ff'$. It follows that
\[
\phi(g)=\phi(f)\phi(f')\leq \phi(f).
\]
This shows a certain order preserving property of $\phi$.

Let $f\in \EX$ be nowhere vanishing. Then there exists a positive
integer $n$ such that $\lambda ^n \leq f$. Since
$\phi(\lambda)=\lambda$ and $\phi$ is multiplicative, it follows
that
\[
\lambda^n=\phi(\lambda^n)\leq \phi(f).
\]
This implies that $\phi(f)$ is also nowhere vanishing. Hence
$\phi$ preserves the nowhere vanishing elements between $E(C(X))$
and $E(C(Y))$. Since $\phi^{-1}$ has similar properties as
$\phi$, we obtain that the above two preserving properties of
$\phi$ hold in both directions.

It is now easy to see that the transformation $\psi: C(X)^+ \to
C(Y)^+$ defined by
\[
\psi(h)= - \ln \phi(\exp(-h)) \qquad (h\in C(X)^+)
\]
is bijective, additive and preserves the order in both
directions. In particular, we infer that $\psi$ is positive
homogeneous. Following the construction in the proof of
Proposition~\ref{P:het1} where we obtain $\Phi$ from $\psi_1$, we
see that our present map $\psi$ extends to a bijective unital
linear transformation $\Phi : C(X) \to C(Y)$ which preserves the
order in both directions. Hence, by the already mentioned result
of Kadison, $\Phi$ is a Jordan *-isomorphism. Because of
commutativity, it means that $\Phi$ is an isomorphism between the
function algebras $C(X)$ and $C(Y)$. The general form of such
transformations is well-known. Namely, we know that there is a
homeomorphism $\varphi :Y \to X$ for which we have
\[
\Phi(h)=h \circ \varphi \qquad (h \in C(X)).
\]
Referring to the relation between $\Phi$ and our original
transformation $\phi$ we easily obtain that $\phi(f)=f \circ
\varphi$ holds for every nowhere vanishing function $f \in \EX$.
We assert that this formula is valid for every element of $\EX$.
In order to see it, consider the transformation
\[
\phi': f \longmapsto \phi(f)\circ \varphi^{-1}.
\]
Clearly, this is a bijective multiplicative selfmap of $\EX$
which acts as the identity on the nowhere vanishing elements of
$\EX$. We prove that the same holds on the whole set $\EX$. In
fact, suppose that $f\in \EX$. It is trivial to see that there is
a sequence $(f_n)$ of nowhere vanishing elements of $\EX$ such
that $f\leq f_n$ and $(f_n)$ converges (uniformly) to $f$. By the
order preserving property of $\phi'$ obtained in the second
paragraph of the proof for $\phi$, we have $\phi'(f)\leq
\phi'(f_n)=f_n.$ Taking limit, we arrive at
\[
\phi'(f)\leq f.
\]
Since ${\phi'}^{-1}$ has similar properties as $\phi'$, we also
have
\[
f={\phi'}^{-1}(\phi'(f))\leq \phi'(f).
\]
These result in $\phi'(f)=f$ $(f\in \EX)$. Therefore, we obtain
\[
\phi(f)=f\circ \varphi \qquad (f \in \EX)
\]
and the proof is complete.
\end{proof}

In the remaining part of the paper $\A, \B$ denote von Neumann
algebras and $\phi: \EA \to \EB$ is a sequential isomorphism.

\begin{lemma}\label{L:het3}
We have $\phi(\PA)=\PB$. The restriction of $\phi$ onto $\PA$ is
a bijective map from $\PA$ onto $\PB$ which preserves the order
and the orthogonality in both directions and hence it is
completely orthoadditive.
\end{lemma}

\begin{proof}
It is trivial that for an arbitrary $A\in \EA$ we have
\[
\text{$A$ is a projection } \Longleftrightarrow A \circ A=A.
\]
This implies that $\phi$ preserves the projections in both
directions.

Let $P,Q\in \PA$. Clearly, we have
\[
P\leq Q \Longleftrightarrow Q\circ P=P.
\]
This implies that $\phi$ preserves the order among projections in
both directions. In particular, we have $\phi(0)=0$ and
$\phi(I)=I$. As for the orthogonality between projections, we
have
\[
P Q=0 \Longleftrightarrow P\circ Q=0.
\]
This implies that $\phi$ preserves the orthogonality between
projections in both directions.

Now, let $(P_\alpha)$ be an arbitrary collection of mutually
orthogonal projections in $\A$ with sum $P$. Using the above
verified properties of $\phi$, we deduce that $(\phi(P_\alpha))$
is a collection of mutually orthogonal projections and we have
$\phi(P_\alpha)\leq \phi(P)$. This implies that
\[
\sum_\alpha \phi(P_\alpha)\leq \phi(P)=\phi(\sum_\alpha
P_\alpha).
\]
Since $\phi^{-1}$ has similar properties as $\phi$, it follows
that
\[
P=\sum_\alpha \phi^{-1}(\phi(P_\alpha))\leq
\phi^{-1}(\sum_\alpha \phi(P_\alpha)).
\]
By the order preserving property of $\phi$ we obtain that
\[
\phi(P)\leq \sum_\alpha \phi(P_\alpha).
\]
Therefore, we have
\[
\sum_\alpha \phi(P_\alpha)= \phi(P)
\]
which means that $\phi$ is completely orthoadditive on $\PA$.
\end{proof}

\begin{lemma}\label{L:het4}
Let $Z,Z' \in \EA$ be central elements in $\A$ and let $P,P'\in
\A$ be mutually orthogonal projections. Then we have
\[
\phi(ZP+Z'P')=\phi(ZP)+\phi(Z'P').
\]
\end{lemma}

\begin{proof}
First observe that $\phi(ZP+Z'P')$ is defined. Indeed, for
$ZP+Z'P'=PZP+P'Z'P'$ we have
\[
0\leq PZP+P'Z'P' \leq PIP+P'IP'=P+P'\leq I,
\]
that is, $ZP+Z'P'\in \EA$. Next, since $ZP+Z'P'$ commutes with
$P,P', P+P'$, the same holds true for its square root. Hence,
using the orthoadditivity of the sequential isomorphism $\phi$ on
$\PA$ (see Lemma~\ref{L:het3}) we can compute
\begin{equation*}
\begin{gathered}
\phi(ZP+Z'P')=\\
\phi((Z P+{Z'} P')^\half (P+P')(Z P+{Z'} P')^\half)=\\
\phi(Z P+{Z'} P')^\half \phi(P+P')\phi(Z P+{Z'} P')^\half=\\
\phi(Z P+{Z'} P')^\half (\phi(P)+\phi(P'))\phi(Z P+{Z'} P')^\half=\\
\phi(Z P+{Z'} P')^\half \phi(P)\phi(Z P+{Z'} P')^\half+
\phi(Z P+{Z'} P')^\half \phi(P')\phi(Z P+{Z'} P')^\half=\\
\phi((Z P+{Z'}P')^\half P(Z P+{Z'} P')^\half)+
\phi((Z P+{Z'}P')^\half P'(Z P+{Z'} P')^\half)=\\
\phi(ZP)+\phi(Z'P').
\end{gathered}
\end{equation*}
The proof is complete.
\end{proof}

\begin{lemma}\label{L:het5}
The sequential isomorphism $\phi$ preserves commutativity in both
directions. Therefore, we have $\phi(E(\mathcal
Z(\A)))=E(\mathcal Z(\B))$.
\end{lemma}

\begin{proof}
The result \cite[Corollary 3]{GudNag02} says that for any Hilbert
space effects $A,B$ we have $A\circ B=B\circ A$ if and only if
$AB=BA$. Since $\phi$ is a sequential isomorphism, this
characterization of commutativity clearly implies that $\phi$
preserves commutativity in both directions. To the second
statement observe that the set of all elements of $\EA$ which
commute with every effect in $\A$ equals to $E(\mathcal Z(\A))$.
\end{proof}

We say that the collection $(E_{ij})$ of operators in the von
Neumann algebra $\A$ forms a self-adjoint system of $n\times n$
matrix units if $n$ is any cardinal number, the index set in
which $i,j$ run has cardinality $n$, for every $i,j,k,l$ we have
\[
E_{ij}E_{kl}=
\left\{%
\begin{array}{ll}
    0, & \hbox{$j\neq k$;} \\
    E_{il}, & \hbox{$j=k$,} \\
\end{array}%
\right.
\]
$\sum_i E_{ii}=I$ in the strong operator topology and
$E_{ij}^*=E_{ji}$ holds for all $i,j$.

\begin{lemma}\label{L:het6}
Suppose that $\A$ has a self-adjoint system of $n\times n$ matrix
units for some cardinal number $n\geq 2$. Then the restriction of
$\phi$ onto $E(\mathcal Z(\A))$ is an $E$-isomorphism onto
$E(\mathcal Z(\B))$. Moreover, $\phi$ is homogeneous in the sense
that
\[
\phi(\lambda A)=\lambda \phi(A)
\]
holds for every $A\in \EA$ and $\lambda \in [0,1]$.
\end{lemma}

\begin{proof}
From Lemma~\ref{L:het5} we know that the restriction of $\phi$
onto $E(\mathcal Z(\A))$ is a sequential isomorphism onto
$E(\mathcal Z(\B))$.

Let $(E_{ij})$ be a self-adjoint system of $n\times n$ matrix
units in $\A$. Let $Z,Z' \in \EA$ be central elements. Suppose
that $Z,Z' \leq \half I$. For temporary use, fix the matrix units
$E=E_{ii}, V=E_{ij}$ $i\neq j$. Define
\[
P=\half (E+V)^*(E+V), \quad
P'=\half (E-V)^*(E-V).
\]
(The trick to use these operators to prove a kind of additivity
of $\phi$ comes from the proof of Lemma 2.1 in \cite{Hak86b}
where the author studied the problem of the additivity of
so-called Jordan *-maps between operator algebras.) It is obvious
that $P,P'$ are mutually orthogonal projections. It follows from
$Z,Z' \leq \half I$ that $Z+Z', 2Z, 2Z'\in \EA$. We have
\begin{equation}\label{E:het3}
E(2ZP)E=2Z(EPE)=ZE, \quad E(2Z'P')E=2Z'(EP'E)=Z'E.
\end{equation}
Since $\phi$ preserves commutativity in both directions,
$\phi(Z), \phi(Z'), \phi(Z+Z')$ are central elements in $\EB$.
Using Lemma~\ref{L:het4} and \eqref{E:het3} we can compute
\[
\begin{gathered}
\phi(Z+Z')\phi(E)=
\phi(Z+Z')^\half \phi(E)\phi(Z+Z')^\half=\\
\phi((Z+Z')E)=
\phi(E(2ZP+2Z'P')E)=\\
\phi(E)\phi(2ZP+2Z'P')\phi(E)=
\phi(E)(\phi(2ZP)+\phi(2Z'P'))\phi(E)=\\
\phi(E)\phi(2ZP)\phi(E)+\phi(E)\phi(2Z'P')\phi(E)=\\
\phi(2EZPE)+\phi(2EZ'P'E)=
\phi(ZE)+\phi(Z'E)=\\
\phi(Z)\phi(E)+\phi(Z')\phi(E)=
(\phi(Z)+\phi(Z'))\phi(E).
\end{gathered}
\]
Therefore, we obtain that
\begin{equation}\label{E:het1}
\phi(Z+Z')\phi(E_{ii})=(\phi(Z)+\phi(Z'))\phi(E_{ii})
\end{equation}
holds for every $i$. We know that $\sum_i E_{ii} =I$ and by the
complete orthoadditivity of $\phi$ (see Lemma~\ref{L:het3}) this
yields that $\sum_i \phi(E_{ii})=I$. Consequently, we deduce from
\eqref{E:het1} that
\[
\phi(Z+Z')=\phi(Z)+\phi(Z')
\]
holds for every central elements $Z,Z'\in \EA$ with $Z,Z' \leq
\half I$. In particular, we obtain that
\[
I=\phi(I)=\phi\biggl(\half I +\half I\biggr)=\phi\biggl(\half
I\biggr) +\phi\biggl(\half I\biggr)
\]
and this implies that
\[
\phi\biggl(\half I\biggr)=\half I.
\]
Restricting $\phi$ onto $E(\mathcal Z (\A))$, we have a
sequential isomorphism between the sets of commutative von
Neumann algebra effects which maps a nontrivial scalar to itself.
Applying Proposition~\ref{P:het2}, it follows that this
restriction of $\phi$ can be extended to a *-isomorphism from
$\mathcal Z(\A)$ onto $\mathcal Z(\B)$. This implies the first
assertion in our statement.

The homogeneity of $\phi$ is now easy to see. In fact, pick
$\lambda \in [0,1]$ and $A\in \EA$. As $\phi$ is homogeneous on
$E(\mathcal Z(\A))$, we compute
\[
\phi(\lambda A)=
\phi(A)^\half \phi(\lambda I) \phi(A)^\half=
\phi(A)^\half \lambda \phi(I) \phi(A)^\half=
\lambda \phi(A).
\]
\end{proof}

\begin{lemma}\label{L:het7}
Suppose that $P\in \A$ is an abelian projection and $A\in \A^+$.
Then there is a positive element $Z$ in $\mathcal Z(\A)$ such
that $PAP=ZP$.
\end{lemma}

\begin{proof}
The operator $PAP$ is a positive element in the $C^*$-algebra
$P\A P$. Therefore, there exists an element $B\in P\A P$ such
that $PAP=B^*B$. But, as $P$ is abelian, we know that
\[
P\A P=\mathcal Z (\A) P
\]
(see \cite[6.4.2. Proposition]{KadRin86}). So, there is a central
element $Z$ in $\A$ such that $B=ZP$. We have
\[
PAP=B^*B=(ZP)^*(ZP)=Z^*Z P
\]
and this verifies our statement.
\end{proof}

\begin{lemma}\label{L:het8}
Let $\A$ be a von Neumann algebra of type I$_n$ with $n <\infty$.
Let $A\in \A$ be such that $PAP=0$ holds for every abelian
projection $P\in \A$. Then we have $A=0$.
\end{lemma}

\begin{proof}
Considering the real and imaginary parts of $A$, we can clearly
assume that $A$ is self-adjoint.

By \cite[8.2.8. Theorem]{KadRin86} there is a center-valued trace
on $\A$ which we denote by $\tau$. In what follows we assume that
$\tau$ has the properties listed in that theorem. We have
\[
0=\tau(PAP)=\tau (APP)=\tau(AP)
\]
for every abelian projection $P$ in $\A$. It is known that any
nonzero projection in a type I algebra is the sum of mutually
orthogonal abelian projections (see \cite[6.4.8. Proposition and
6.5.1. Definition]{KadRin86}). Hence, if we pick an arbitrary
projection $Q$ in $\A$, then there is a collection $(P_\alpha)$
of mutually orthogonal abelian projections such that
$Q=\sum_\alpha P_\alpha$. As the weak and ultraweak topologies
coincide on the unit ball of any von Neumann algebra (see the
discussion in \cite[7.4.4. Remark]{KadRin86}), it follows that
$AQ=\sum_\alpha AP_\alpha$ holds in the ultraweak topology. By
the ultraweak continuity of $\tau$ we infer that
\begin{equation}\label{E:het6}
\tau(AQ)=\sum_\alpha \tau(AP_\alpha)=0.
\end{equation}
Since the linear span of the set of all projections in a von
Neumann algebra is dense with respect to the norm topology and
$\tau$ is norm continuous, it follows from \eqref{E:het6} that
$\tau(A^2)=0$. By the definiteness of the trace, we obtain that
$A^2=0$. This implies $A=0$ and the proof is complete.
\end{proof}

\begin{lemma}\label{L:het9}
The sequential isomorphism $\phi$ preserves the abelian
projections and the equivalence among them in both directions.
\end{lemma}

\begin{proof}
Clearly, $P\in \A$ is abelian if and only if the set $P\EA P$ is
commutative. As $\phi$ is a sequential isomorphism, it preserves
commutativity in both directions (see Lemma~\ref{L:het5}).
Therefore, the set $P\EA P$ is commutative if and only if so is
its image under $\phi$. But this image equals
\[
\phi(P\EA P)=\phi(P)\phi(\EA)\phi(P)=\phi(P)\EB \phi(P).
\]
Therefore, $P$ is abelian if and only if $\phi(P)$ is abelian.

As for the preservation of the equivalence between abelian
projections, we recall that two abelian projections are
equivalent if and only if their central carriers coincide (see
\cite[6.2.8. Proposition and 6.4.6 Proposition]{KadRin86}). But
the notion of the central carrier is expressed by order and
commutativity both of them being preserved by $\phi$ in both
directions. This implies the second assertion of our statement.
\end{proof}

After this preparation, now we are in a position to prove that
every sequential isomorphism between the sets of effects in type
I$_n$ algebras $(2\leq n<\infty)$ is an $E$-isomorphism.

\begin{proposition}\label{P:het10}
Suppose that $\A, \B$ are type I$_n$ algebras with $2\leq
n<\infty$. Then $\phi$ is an $E$-isomorphism.
\end{proposition}

\begin{proof}
Let $A,A'\in \EA$ be such that $A+A' \in \EA$. Pick an arbitrary
abelian projection $P\in \A$. By Lemma~\ref{L:het7} we have
positive central elements $Z,Z'$ in $\A$ such that
\begin{equation}\label{E:het2}
PAP=ZP, \quad PA'P=Z'P.
\end{equation}
Clearly, we can choose a scalar $\lambda \in ]0,1]$ such that
$0\leq \lambda (Z+Z')\leq I$. It is well-known that for $n\geq
2$, any type I$_n$ algebra has a self-adjoint system of $n\times
n$ matrix units (see \cite[6.6.3. Lemma]{KadRin86}). Therefore,
Lemma~\ref{L:het6} applies and using \eqref{E:het2} we can
compute
\[
\begin{gathered}
\lambda \phi(P)\phi(A+A')\phi(P)=\\
\lambda \phi(P(A+A')P)=
\phi(\lambda P(A+A')P)=
\phi(\lambda (ZP+Z'P))=\\
\phi(P) \phi(\lambda (Z+Z'))\phi(P)=
\phi(P) (\phi(\lambda Z)+\phi(\lambda Z'))\phi(P)=\\
\phi(P) \phi(\lambda Z)\phi(P) +\phi(P) \phi(\lambda Z')\phi(P)=
\phi(\lambda PZP) + \phi(\lambda PZ'P)=\\
\phi(\lambda PAP) + \phi(\lambda PA'P)=
\lambda(\phi(PAP) + \phi(PA'P))=\\
\lambda \phi(P)(\phi(A) + \phi(A'))\phi(P).
\end{gathered}
\]
Consequently, we obtain that
\[
\phi(P)\phi(A+A')\phi(P)=\phi(P)(\phi(A) + \phi(A'))\phi(P)
\]
holds for every abelian projection $P\in \A$. Referring to
Lemma~\ref{L:het9} we infer that
\[
Q(\phi(A+A')-(\phi(A) + \phi(A'))Q=0
\]
is valid for every abelian projection $Q$ in $\B$. By
Lemma~\ref{L:het8} this implies that
\[
\phi(A+A')=\phi(A) + \phi(A').
\]
In particular, we also get that $\phi(A)+\phi(A')\in \EB$. Since
$\phi^{-1}$ has similar properties as $\phi$, we obtain that
$\phi$ is an $E$-isomorphism.
\end{proof}

In the case of von Neumann algebras of the remaining types we can
follow a different approach based on the solution of the
so-called Mackey-Gleason problem due to Bunce and Wright. We have
the following result.

\begin{proposition}\label{P:het11}
Suppose that $\A, \B$ are both of type I$_n$ with $2<n$ or $\A,
\B$ have no type I direct summands. Then $\phi$ is an
$E$-isomorphism.
\end{proposition}

\begin{proof}
We know from Lemma~\ref{L:het3} that $\phi$, when restricted onto
the set of all projections in $\A$, is orthoadditive. We now
apply the deep result of Bunce and Wright \cite{BunWri92} stating
that every bounded orthoadditive map from the set of all
projections of a von Neumann algebra without type I$_2$ direct
summand into a Banach space can be extended to a continuous
linear transformation defined on the whole algebra. Therefore, we
have a continuous linear transformation $L: \A \to \B$ such that
\[
L(P)=\phi(P) \qquad (P\in P(\A)).
\]
We claim that $L$ coincides with $\phi$ on the whole set $\EA$.
To see this, first observe that
\[
L(\lambda P)=\phi(\lambda P)
\]
whenever $\lambda \in [0,1], P\in P(\A)$. Indeed, this follows
from the homogeneity of $\phi$ which was asserted in
Lemma~\ref{L:het6} (notice that by \cite[6.5.6. Lemma and 6.6.4.
Lemma]{KadRin86}, $\A$ has a self-adjoint system of $n\times n$
matrix units). Now, let $(P_i)$ be a finite collection of
mutually orthogonal projections in $\A$ with sum $I$ and pick
scalars $\lambda_i \in [0,1]$. Since $\sum_i \phi(P_i)=I$, we can
compute
\[
\begin{gathered}
\phi(\sum_i \lambda_i P_i)=
\phi(\sum_i \lambda_i P_i)^\half \sum_k \phi(P_k)\phi(\sum_i \lambda_i P_i)^\half=\\
\sum_k \phi(\sum_i \lambda_i P_i)^\half \phi(P_k)\phi(\sum_i \lambda_i P_i)^\half=
\sum_k \phi((\sum_i \lambda_i P_i)^\half P_k (\sum_i \lambda_i P_i)^\half)=\\
\sum_k \phi(\lambda_k P_k)= \sum_k L(\lambda_k P_k)= L(\sum_k \lambda_k P_k).
\end{gathered}
\]
Let $A\in \EA$ be arbitrary. It follows easily from the spectral
theorem of normal operators and the properties of the spectral
integral that for any $\epsilon
>0$ there are operators $A_l, A^u \in \EA$ of the form $\sum_i
\lambda_i P_i$ (where $\lambda_i, P_i$ are like above) and
operators $A_l', {A^u}'\in \EA$ such that the whole set $\{A,
A_l, A^u, A_l', {A^u}'\}$ is commutative, $\| A_l -A\|, \| A^u
-A\| <\epsilon$, and
\[
A_l=AA_l', \quad A=A_u{A_u}'.
\]
We compute
\[
L(A_l)=
\phi(A_l)=
\phi(A^\half A_l' A^\half)=
\phi(A)^\half \phi(A_l') \phi(A)^\half\leq
\phi(A)^\half I \phi(A)^\half=
\phi(A)
\]
and one can prove in a similar manner that
\[
\phi(A)\leq \phi(A_u)=L(A_u).
\]
Therefore, we have $L(A_l)\leq \phi(A)\leq L(A^u)$ and by the
continuity of $L$ we infer that
\[
L(A)\leq \phi(A) \leq L(A).
\]
Consequently, we have $\phi(A)=L(A)$ for any $A\in \EA$.

It follows that whenever $A,A' \in \EA$ are such that
$A+A'\in\EA$, we have
\[
\phi(A+A')=L(A+A')=L(A)+L(A')=\phi(A)+\phi(A').
\]
In particular, we have $\phi(A)+\phi(A')\in \EB$. Just as in the
proof of Proposition~\ref{P:het10}, referring to the fact that
$\phi^{-1}$ has similar properties as $\phi$, we obtain that
$\phi$ is an $E$-isomorphism.
\end{proof}

Putting together all the information that we have collected so
far, it is now an easy task to prove our main result.

\begin{proof}[Proof of Theorem~\ref{T:hetfo}]
Let $P\in \A$ be a nonzero projection. Then $\phi(P)$ is also a
nonzero projection and the restriction of $\phi$ onto $P\EA
P=E(P\A P)$ gives rise to a sequential isomorphism from $E(P\A
P)$ onto $E(\phi(P)\B \phi(P))$. By the properties of the
sequential isomorphisms formulated in Lemmas~\ref{L:het3},
\ref{L:het5}, \ref{L:het9}, it follows that $P\A P$ is of type I,
or of type I$_n$, or has no direct summand of type I, or has no
commutative direct summand if and only if the same holds for
$\phi(P)\B \phi(P)$.

Consider the type decomposition of $\A$ (see, for example,
\cite[6.5.2. Theorem]{KadRin86}). For any cardinal number $n$
(not exceeding the dimension of the Hilbert space on which the
elements of $\A$ act) let $P_n\in \A$ be a central projection
such that the algebra $\A P_n$ is of type $I_n$ (or $P_n=0$) and
in case $Q=I-\sum_n P_n \neq 0$ the algebra $\A Q$ has no type I
direct summand. It follows from the first paragraph of the proof
that the collection $(\phi(P_n))$ of central projections in $\B$
has the same properties (relating the algebra $\B$ in the place
of $\A$, of course).

Clearly, the algebras $\A_1=\A P_1$ and $\B_1=\B \phi(P_1)$ are
commutative and their direct complements in $\A$, respectively in
$\B$ have no commutative direct summands.

Apply our auxiliary result Proposition~\ref{P:het10} for the
direct summands $\A P_n$ and $\B \phi(P_n)$ whenever $2\leq
n<\infty$. Moreover, apply Proposition~\ref{P:het11} for the
direct summands $\A P_n$ and $\B \phi(P_n)$ whenever $n$ is an
infinite cardinal and do the same for $\A Q$ and $\B \phi(Q)$.
Taking direct sums, we obtain that the restriction of $\phi$ onto
$\EA (I-P_1)=E(\A (I-P_1))$ is an $E$-isomorphism onto $E(\B
(I-\phi(P_1)))$. By Proposition~\ref{P:het1} this can be extended
to a Jordan *-isomorphism from $\A(I-P_1)$ onto $\B
(I-\phi(P_1))$. It is well-known that every Jordan *-isomorphism
between two von Neumann algebras induces a direct sum
decomposition of the underlying algebras according to which the
Jordan *-isomorphism under consideration is the direct sum of a
*-isomorphism and a *-antiisomorphism (see \cite[10.5.26.
Exercise]{KadRin86}). This verifies the existence of a
decomposition having the properties (i)-(ii) in our theorem.

As for the last assertion, observe that if $P_1\neq 0$ (i.e.,
when $\A_1, \B_1$ 'really appear') and $\phi$ maps a nontrivial
scalar to itself, then by Proposition~\ref{P:het2}, the
sequential isomorphism $\phi_1$ from $E(\mathcal A_1)$ onto
$E(\mathcal B_1)$ can also be extended to a *-isomorphism from
$\A_1$ onto $\B_1$. The proof is complete.
\end{proof}

\bibliographystyle{amsplain}

\begin{thebibliography}{99}

\bibitem{BunWri92}
L.J. Bunce and J.D.M. Wright,
\emph{The Mackey-Gleason problem},
Bull. Amer. Math. Soc. \textbf{26} (1992), 288--293.

\bibitem{BusLahMit91}
P. Busch, P.J. Lahti and P. Mittelstaedt,
\emph{The Quantum Theory of Measurement,}
Springer-Verlag, 1991.

\bibitem{CasDeVLahLev}
G. Cassinelli, E. De Vito, P. Lahti and A. Levrero,
\emph{Symmetry groups in quantum mechanics and the theorem of
Wigner on the symmetry transformations,} Rev. Math. Phys.
\textbf{8} (1997), 921--941.

\bibitem{GudGre02}
S. Gudder and R. Greechie,
\emph{Sequential products on effect algebras,}
Rep. Math. Phys. \textbf{49} (2002), 87--111.

\bibitem{GudNag01}
S. Gudder and G. Nagy,
\emph{Sequential quantum measurements,}
J. Math. Phys. \textbf{42} (2001), 5212--5222.

\bibitem{GudNag02}
S. Gudder and G. Nagy,
\emph{Sequentially independent effects,}
Proc. Amer. Math. Soc. \textbf{130} (2002), 1125--1130.

\bibitem{Hak86b}
J. Hakeda,
\emph{Additivity of Jordan *-maps on $AW^*$-algebras},
Proc. Amer. Math. Soc. \textbf{96} (1986), 413--420.

\bibitem{KadRin86}
R.V. Kadison and J.R. Ringrose,
\emph{Fundamentals of the Theory of Operator Algebras, Vol II.,}
Academic Press, 1986.

\bibitem{Kra83}
K. Kraus,
\emph{States, Effects and Operations,}
Lecture Notes in Physics, Vol. 190, Springer-Verlag, 1983.

\bibitem{Lud83}
G. Ludwig,
\emph{Foundation of Quantum Mechanics I,}
Springer-Verlag, 1983.

\bibitem{ML00f}
L. Moln\'ar,
\emph{On some automorphisms of the set of effects on Hilbert space},
Lett. Math. Phys. \textbf{51} (2000), 37--45.

\bibitem{MLxx7}
L. Moln\'ar,
\emph{Preservers on Hilbert space effects,}
Linear Algebra Appl., to appear.

\bibitem{MLxx5}
L. Moln\'ar and P. \v Semrl, \emph{Conditional affine and
conditional sequential automorphisms of the set of Hilbert space
effects,} preprint.

\end{thebibliography}

\end{document}